\documentclass{amsart}
\usepackage{palatino,hhline, bbm}
\usepackage{arydshln} 
\usepackage{amsthm, amsmath}
\usepackage{amssymb} 
\usepackage{amscd}
\usepackage{mathrsfs}
\renewcommand{\mathcal}{\mathscr}
\usepackage[hmargin=3cm, vmargin=3cm]{geometry}
\usepackage[breaklinks,bookmarksopen,bookmarksnumbered]{hyperref}
\usepackage[all]{xy}

\theoremstyle{plain}

\newtheorem*{lem}{Lemma}
\newtheorem{prop}{Proposition}
\newtheorem*{cor}{Corollary}

\theoremstyle{remark}

\newcommand\pr{\noindent\textit{Proof} : }

\newcommand\rond{\kern 1pt{\scriptstyle\circ}\kern 1pt}

\def\lr#1{\langle {#1} \rangle}

\newcommand\Ker{\operatorname{Ker}}

\newcommand{\mo}{\smallsetminus}

\newcommand\Z{\mathbb{Z}}

\newcommand\F{\mathbb{F}}

\newcommand\OO{\operatorname{O} }

\newcommand\iso{\vbox{\hbox to .8cm{\hfill{$\scriptstyle\sim$}\hfill}
\nointerlineskip\hbox to .8cm{{\hfill$\longrightarrow $\hfill}} }}
\newcommand\bir{\vbox{\hbox to .8cm{\hfill{$\scriptstyle\sim$}\hfill}
\nointerlineskip\hbox to .8cm{{\hfill$\dasharrow $\hfill}} }}

\begin{document}
\title{Reduction  mod. 2 of del Pezzo lattices}
\author[Arnaud Beauville]{Arnaud Beauville}
\address{Universit\'e C\^ote d'Azur\\
CNRS -- Laboratoire J.-A. Dieudonn\'e\\
%UMR 7351 du CNRS\\
Parc Valrose\\
F-06108 Nice cedex 2, France}
\email{arnaud.beauville@unice.fr}
 
\begin{abstract}
For an even lattice $L$, the form $v\mapsto \frac{1}{2}v^2 $ induces a quadratic form $q$ on the $\F_2$-vector space $L/2L$. For the lattices associated to some particular root systems, we show that reduction mod.\ 2 induces a bijection  between the roots of $L$ and the vectors of $L/2L$ with $q=1$, and an isomorphism  $\OO(L)/\{\pm 1\}\iso \OO(L/2L)$.
\end{abstract}
%\subjclass[2010]{Primary: 14M20; Secondary: 14E08, 14J45}
\maketitle 

\section{Introduction}

Let  $L$ be an even lattice. The form $v\mapsto \frac{1}{2}v^2 $ on $L$ induces a quadratic form $q$ on $L_2:=L/2L$; we want to investigate what information on $L$ we get from the quadratic $\F_2$-vector space $L_2$. In general not much can be said, but for the lattices associated to certain root systems (specifically $A_1\times A_2, A_4,D_5,E_6,E_7,E_8$),
 the relation is tight: the  roots correspond to the vectors of $L_2$ with $q=1$ (with a grain of salt for $E_7$), and the group $\OO(L)$ is very close to the orthogonal group $\operatorname{O}(L_2) $. We call the corresponding lattices \emph{del Pezzo lattices}, since they appear naturally in the Picard group of del Pezzo surfaces \cite[\S\,25]{M}, but this will play no role in this note.

\medskip	
$\bullet$ Added in version 2: Igor Dolgachev pointed out the reference \cite{G}, which contains much more complete results (with a much more elaborate technology). I leave this preprint as a possible elementary introduction to the subject.

\medskip	
\section{Reduction mod. $2$ of roots}

Let $L$ be an even lattice.
We endow $L_2:=L/2L$ with the bilinear form $(\ |\ )$ deduced from the scalar product on $L$, and the quadratic form $q$ deduced from the
 form  $v\mapsto \frac{1}{2}v^2 $. We have $q(x+y)=q(x)+q(y)+(x|y)$ for $x,y$ in $L_2$.
 
 \begin{lem}
 Let $R$ be a root system of type $A, D,E$, and $L$ the corresponding root lattice.  

$1)$ The inclusion $R\subset L$ induces an injection $R/\{\pm 1\} \hookrightarrow L_2$. 

$2)$ The reduction map $\rho :\OO(L)\rightarrow \OO(L_2)$ induces an injection $\OO(L)/\{\pm 1\}\hookrightarrow \OO(L_2) $.

\end{lem}

\pr 1) Let $\alpha ,\beta \in R$ such that $\alpha -\beta =2\lambda $, with $\lambda \in L$, $\lambda \neq 0$. Then $8\leq (\alpha -\beta )^2=4-2(\alpha \cdot \beta )$, hence $(\alpha \cdot \beta )\leq -2$; by the Cauchy-Schwarz inequality this implies $\beta =-\alpha $.

\smallskip	
2) Let $u\in \OO(L)$ such that $\rho (u)=1$. For each root $\alpha $, $u(\alpha )$ is a root, with $u(\alpha )-\alpha \in 2L$. By 1), this implies $u(\alpha )=\varepsilon  \alpha $, with $\varepsilon =\pm 1$. Let $\beta $ be another root with $(\alpha \cdot \beta) \neq 0$; since $(u(\alpha )\cdot u(\beta ))=(\alpha \cdot \beta) $, we have $u(\beta  )=\varepsilon \beta $. Taking for $\alpha $ a simple root, it follows that $u(\gamma )=\varepsilon \gamma $ for all simple roots $\gamma $, hence $u=\varepsilon \cdot 1_L$.\qed

\medskip	
From now on, we will only consider a particular type of lattices, which we call \emph{del Pezzo lattices}.  For $3\leq n\leq 8$, let $L_0$ be the free $\Z$-module  with basis $E_0,\ldots ,E_n$, endowed with the Lorentzian scalar product defined by $E_0^2=-1$, $E_i^2=1$ for $i\geq 1$, and $(E_i\cdot E_j)=0$ for $i\neq j$. Let $K$ be the element $3E_0-\sum\limits_{i\geq 1}E_i$. We put $L=K^{\perp}\subset L_0$; this is an even lattice, of discriminant $9-n$. The vectors of square 2 form a root system $R\subset L$, of type $A_1\times A_2, A_4,D_5,E_6,E_7,E_8$ for $n=3,\ldots ,8$ \cite[\S\,25]{M}. 

\smallskip	
Let $\varphi :R\rightarrow L_2$ be the reduction mod. 2. By  the Lemma, $\varphi $ induces an injective map $R/\{\pm 1\}\hookrightarrow L_2 $; its image is contained in $q^{-1}(1)$.

\begin{prop}\label{root}
$\varphi$ induces a bijection  from $R/\{\pm 1\}$ onto $q^{-1}(1) $ if $n\neq 7$, and onto  $q^{-1}(1)\mo\{k\} $   if $n=7$, where $k$ is the nonzero element of the kernel of $(\ |\ )$.
\end{prop}

\pr Let $e_0,\ldots ,e_n$ be the images of $E_0,\ldots ,E_n$ in $L_0/2L_0$; then $L_2$ is the orthogonal in $L_0/2L_0$ of the vector $k=\sum e_i$. We have
\[q(e_0+e_i)=0\, ,\quad q(e_i+e_j)=1\quad \mbox{for }0<i<j\,.\]The   bilinear form $(\ |\ )$ on $L_2$ is non-degenerate if $n$ is even; if $n$ is odd its kernel is $\{0,k\} $. 

\smallskip

Let $v\in q^{-1}(1)$. Then $v=\sum\limits_{i\in I}e_i$, with $\#I=m$ even. We put $E_I:=\sum\limits_{i\in I}E_i$. 

a)  Suppose $0\notin I$. If $m=2$, $v=e_i+e_j=\varphi (E_i-E_j)$. If $m=6$, $v=\varphi (2E_0-E_I)$. 

b) Suppose $0\in I$. If $m=4$, $v=\varphi (2E_0-E_I)$. If $m=n=8$, $v=\varphi (4E_0-E_I -2E_\ell)$, with $\ell\notin I$.  

c) It remains to show that for $n=7$,  
the vector $k=\sum e_i$  cannot be written $\varphi (\alpha )$ for some $\alpha \in R$. Indeed we would have $\alpha =\sum n_iE_i$ with each $n_i$ odd; but then $n_i^2\equiv 1\ \operatorname{mod.}\,8 $, hence $\alpha ^2\equiv 6\ \operatorname{mod.}\,8$, contradicting $\alpha ^2=2$.\qed

 \bigskip	
 It is instructive to check the equality $\frac{1}{2} \#R=\#q^{-1}(1)$ (or $\#q^{-1}(1)-1$ for $n=7$) in each case:
 
 \smallskip	
\noindent a) $n=2m$.
 
 In this case the form $(\ |\ )$ is non-degenerate; the order of $q^{-1}(1)$ is $2^{m-1}(2^m +(-1)^{a})$, where $a$ is the \emph{Arf invariant}  $\operatorname{Arf}( q)$: if $(\delta _1,\ldots ,\delta _m; \varepsilon _1,\ldots ,\varepsilon _m)$ is a symplectic basis of $L_2$, $a=\sum q(\delta _i)q(\varepsilon _i)$. 
  
We can take for instance  $\delta _k:= \sum\limits_{i=0}^{2k-1}e_i$, $\varepsilon _k=\sum\limits_{i=0}^{2k-2}e_i+e_{2k} $. Then $q(\delta _k)=q(\varepsilon _k)= 0$ for $k$ odd, $1$ for $k$ even, hence $\operatorname{Arf}(q)= 0$ for $n=4$ or $6$, $\operatorname{Arf}(q)= 1$ for $n=8$. Therefore we get $\#q^{-1}(1)=10, 36, 120=\frac{1}{2}\# R $ for $n=4,6,8$.

\medskip	
\noindent b) $n=3$ or $7$.

For $n$ odd, the vector $k$ is in $L_2$, and spans the kernel of the form $(\ |\ )$.  For $n=3$ or $7$ we have $q(k)=1$. The involution $v\mapsto v+k $ exchanges $q^{-1}(1)$ and $q^{-1}(0)$, hence $\#q^{-1}(1)=\frac{1}{2}\# L_2= 2^{n-1} $. Indeed $\#(A_1\times A_2)= 2+6=8$, and $\# E_7= 126= 2(2^6-1)$.  

\medskip	
\noindent c) $n=5$.

Here $q(k)=0$. The form $q$ descends to a non-degenerate form $\bar{q}$ on $L_2/\lr{k}$; we claim that $L_2/\lr{k}$ is isomorphic to the quadratic $\F_2$-space $L'_2$ associated to $A_4$. Indeed, the subspace $N$ of $L_2$ spanned by $e_i+e_j$ for $i<j<5$, with the quadratic form induced by $q$, is isomorphic to $L'_2$, and the quotient map $N\rightarrow L_2/\lr{k}$ is an isomorphism.
Thus  $\#q^{-1}(1)=2\#\bar{q}^{-1}(1)=\#A_4=20=\frac{1}{2}\# D_5$.

\medskip	
\section{Orthogonal groups}
\begin{prop}
Let $L$ be a del Pezzo lattice of rank $\geq 4$. The reduction map $\rho :\operatorname{O} (L)\rightarrow \operatorname{O}(L_2) $ induces an isomorphism of $\operatorname{O}(L)/\{\pm 1\}  $ onto $\operatorname{O}(L_2) $.
\end{prop}

\pr The induced map is injective by the Lemma, so we have to prove that $\rho $ is surjective.

\smallskip	
\noindent a) $n=4,6,8$.

In this case  the form $(\ |\ )$ is non-degenerate. By \cite[\S\,10]{D},  the group $\OO(L_2)$ is generated by the reflections $r_v: x\mapsto x+(x|v)v$ for $v$ in $L_2$ with $q(v)=1$ -- except for one case that we will consider below. By Proposition 1, we have $v=r(\alpha )$ for some root $\alpha $; then $r_v$ is the reduction mod. 2 of the reflection with respect to $\alpha $. 

The exception  occurs for $n=4$, if   $L_2$ contains a 2-dimensional subspace $N$ such that $q_{|N}=0$ (\emph{loc. cit.}). This does not happen in our case, since the nonzero vectors $v$ with $q(v)=0$ are $f_0:=k+e_0$ and $f_i:=e_0+e_i$ for $1\leq i\leq 4$, and they satisfy $(f_i|f_j)=1$ for $i\neq j$. \qed

\medskip	
\noindent b) $n=3,7$. 

In this case the nonzero element $k$ of $\Ker(\ |\ )$ satisfies $q(k)=1$. Let $H$ be a hyperplane in $L_2$ such that $L_2=H\oplus \F_2k$, and let $p_H$ be the projector from $L_2$ onto $H$. The restriction of $(\ |\ )$ to $H$ is a symplectic form, and  the map $u\mapsto p_H\circ u_{|H}$ is an isomorphism  from $\OO(L_2)$   onto $\operatorname{Sp}(H) $ \cite[\S\,7]{D}; a transvection $x\mapsto (x|v)v$ in $\operatorname{Sp}(H) $ corresponds to the reflection with respect to $v+(1+q(v))k$. Since $\operatorname{Sp}(H) $ is generated by transvections, $\OO(L_2)$ is generated by reflections, and we conclude as in a) (note that  the reflection with respect to $k$ is the identity). 

\medskip	
\noindent c) $n=5$.

Here it is simpler to use a counting argument. The group $\OO(L)/\{\pm 1\} $ is isomorphic to the Weyl group of $D_5$, that is, the semi-direct product $(\Z/2)^4\rtimes \mathfrak{S}_5$ \cite{B}. As we have seen in c) above, the space $L_2/\{k\} $ with the form induced by $q$ is isomorphic to the quadratic $\F_2$-space $L'_2$ for $n=4$; the quotient map induces a surjective map $\pi :\OO(L_2)\rightarrow \OO(L'_2)$. The kernel of $\pi $ consists of all automorphisms of $L_2$ of the form $x\mapsto x+\ell(x)k$, with $\ell\in L_2^*$, $\ell(k)=0$;  it is isomorphic to $(\Z/2)^4$. By a), the group $\OO(L'_2)$ is isomorphic to the Weyl group $\mathfrak{S}_5$ of $A_4$. The result follows.\qed

\begin{cor}
$\rho $  induces an isomorphism of the Weyl group $W$ of $L$ onto $\OO(L_2)$
for $n=4,5 $ or $6$, and of $W/\{\pm 1\} $ onto $\OO(L_2)$ for $n=7$ or $8$.
\end{cor}
This follows from the relation between $W$ and $\OO(L)$, see for instance \cite{B}.

\medskip	
\noindent\emph{Remarks}$.-\ 1)$ For $n=3$, $\rho $ is surjective, with kernel $\{\pm 1\}\times \{\pm 1\}  $. More precisely,
we have $L=L'\oplus L''$, where $L'$ and $L''$ are lattices of type $A_1$ and $A_2$; then $\OO(L)=\OO(L')\times \OO(L'')\allowbreak =\{\pm 1\}\times (\mathfrak{S}_3\times \{\pm 1\} ) $, and $\OO(L_2)= \OO(L_2')\times \OO(L_2'')= \{1\}\times \mathfrak{S}_3 $. 

\smallskip	
2) Propositions 1 and 2 do not hold for other root systems. For $A_n$ for instance, the number of roots is $n(n+1)$, while $\#q^{-1}(1)$ is roughly $2^{n-1}$. Similarly the order of $\OO(L_2)$ is roughly $2^{n^2/2}$, much larger than $n!$.


\bigskip	
\begin{thebibliography}{X}

\bibitem[B]{B} N. Bourbaki\,: \textsl{Lie groups and Lie algebras}, Ch.\ 6. Springer-Verlag, Berlin, 2002.

\bibitem[D]{D} J. Dieudonn\'e\,: \textsl{La g\'eom\'etrie des groupes classiques}.
 Ergebnisse der Math. und ihrer Grenzgebiete \textbf{5}, 3rd ed. Springer-Verlag, Berlin-G\"ottingen-Heidelberg, 1971.
 
 \bibitem[G]{G} R. Griess\,: \textsl{Quotients of Infinite Reflection Groups}.
 Math. Ann. \textbf{263} , 267-278 (1983).
 
 \bibitem[M]{M} Y. Manin\,: \textsl{Cubic forms: Algebra, geometry, arithmetic}.  North-Holland Math. Library \textbf{4}, 2nd ed. North-Holland, Amsterdam, 1986.
\end{thebibliography}
\end{document}